\documentclass{article} 
\usepackage{iclr2024_conference,times}

\usepackage{amsmath,amsfonts,bm}

\def\eqref#1{equation~\ref{#1}}

\def\1{\bm{1}}

\DeclareMathAlphabet{\mathsfit}{\encodingdefault}{\sfdefault}{m}{sl}
\SetMathAlphabet{\mathsfit}{bold}{\encodingdefault}{\sfdefault}{bx}{n}

\newcommand{\E}{\mathbb{E}}

\usepackage{hyperref}
\usepackage{url}
\usepackage{amsthm}
\newtheorem{theorem}{Theorem}[section]

\usepackage{amsmath}

\title{Group Action Approaches in Erdos Quotient Set Problem}

\author{Will Burstein
}

\iclrfinalcopy 

\begin{document}

\maketitle

\begin{abstract}
Let $\mathbb{F}_q$ denote the finite field of $q$ elements.  For $E \subset \mathbb{F}_q^d$, denote the distance set $\Delta(E)= \{\|x-y\|^2:=(x_1-y_1)^2+ \cdots + (x_d-y_d)^2 : (x,y)\in E^2 \}$.
    The Erdos quotient set problem was introduced in \cite{Iosevich_2019} where it was shown that for even $d\geq2$ that if $|E| \subset \mathbb{F}_q^2$ such that $|E| >> q^{d/2}$, then $\frac{\Delta(E)}{\Delta(E)}:= \{\frac{s}{t}:s,t \in \Delta(E), t\not=0\} =\mathbb{F}_q^d$.  The proof of the latter result is quite sophisticated and in \cite{pham2023group}, a simple proof using a group-action approach was obtained for the case of $q \equiv 3 \mod 4$ when $d=2$.  In the $q \equiv 3 \mod 4$ setting, for each $r \in (\mathbb{F}_q)^2$, \cite{pham2023group} showed if $E \subset \mathbb{F}_q$, then $V(r):= \#  \left\{ (a,b,c,d) \in E^2: \frac{\|a-b\|^2}{\|c-d\|^2} = r \right\} >> \frac{|E|^4}{q}$.  In this work we use group action techniques in the $q \equiv 3 \mod 4$ setting, for $d=2$ and improve the results of \cite{pham2023group} by removing the assumption on $r \in (\mathbb{F}_q)^2$.  Specifically we show if $d=2$ and $q \equiv 3 \mod 4$, then for each $r \in \mathbb{F}_q^*$,
$V(r)\geq \frac{|E|^4}{2q}$
if $|E|\geq \sqrt{2}q$
for all $r \in \mathbb{F}_q$.  Finally, we improve the main result of \cite{bhowmik2023near} using our proof techniques from our quotient set results.
\end{abstract}

\section{Introduction}
Let $q=p^n$ for some prime $p$ and $n\geq1$.  $\mathbb{F}_q$ denotes the finite field of $q$ elements.  Given $E \subset \mathbb{F}_q^d$, the Erdos Distance problem asks what that best lower bound of 
\begin{equation*}
\Delta(E):=
\left\{\|x-y \|^2\ = 
(x_1-y_1)^2+ \cdots
+(x_d-y_d)^2:x,y \in E
\right\}
\end{equation*}
is.  The work of \cite{Iosevich2005} formulated the problem and showed that if $|E|\geq 2q^{(d+1)/2}$, then $\Delta(E) = \mathbb{F}_q$.  \cite{hart2007} showed that the $(d+1)/2$ exponent is sharp when $d$ is odd.  Moreover, if $d$ is even, it is shown in \cite{hart2007} that the exponent must be at least $d/2$, for the possibility of $|\Delta(E)| >> q$.  It is conjectured that if $d$ is even and $|E| >> q^{d/2}$, then  $|\Delta(E)| >> q$.  In the more general $k$-simplex setting, \cite{bennett2013group} the results of \cite{Iosevich2005} are generalized bringing in group-action machinery.

There have been various special cases where the $(d+1)/2$ exponent has been improved.  The current results are in the case when $d=2$. \cite{Chapman2011PinnedDS} proved that if $q \equiv 3 \mod 4$, $E \subset \mathbb{F}_q^2$, and $|E| >> q^{4/3}$, then
$|E| > q/2$.  In the more general framework of \cite{bennett2013group}, the latter result was achieved with the removal of the $q \equiv 3 \mod 4$ assumption.  In \cite{Rudnev5_4}, when $q=p$, $p$ an odd prime, the exponent has been improved to $5/4$.  

A similar problem to the Erdos distance problem is the Erdos quotient problem.  The Erdos quotient set is defined as
\begin{equation*}
\frac{\Delta(E)}{\Delta(E)}:=
\left\{\frac{a}{b}: a \in \Delta(E), b\in \Delta(E)^*
\right\}.
\end{equation*}
\cite{Iosevich_2019} formulated the Erdos quotient problem and obtained the result below using Fourier Analytic techniques.
\begin{theorem} 
\label{quotientset}
(Theorem 1.1 of \cite{Iosevich_2019})
Let $E \subset \mathbb{F}_q^d$ $d$ even.  Then if $|E|\geq 9q^{\frac{d}{2}}$, we have 
\begin{equation*}
\mathbb{F}_q = \frac{\Delta(E)}{\Delta(E)}.
\end{equation*}
If $d\geq 3$ is odd and $|E| \geq 6 q^{q/2}$, then
\begin{equation*}
(\mathbb{F}_q^d)^2 \subset \frac{\Delta(E)}{\Delta(E)}.
\end{equation*}
\end{theorem}
\cite{Iosevich_2019} prove both their results of Theorem \ref{quotientset} are sharp in the exponent, $d$. 

One drawback of the result of Theorem \ref{quotientset} is that the proof is quite sophisticated and for a fixed $r$, does not say how many distinct $(a,b,c,d) \in E^4$ such that $\|a-b\|^2/\|c-d\|^2=r$.  In a restricted setting, the next result below alleviates the latter two drawbacks with a simple proof.  The proof doesn't use Fourier analysis and instead uses a group-action approach.

\begin{theorem}
\label{quotient}
(Theorem 1.2 of \cite{pham2023group})  Let $E \subset \mathbb{F}_q^2$ with $q \equiv 3 \mod 4$.  Assume that
$|E| >> q$.  Then for each $r \in (\mathbb{F}_q)^2$, the number of quadruples $(a,b,c,d) \in E^4$ such that 
$\|a-b\|^2/\|c-d\|^2=r$ is at least 
\begin{equation*}
>> |E|^4/q.
\end{equation*}
In particular,
\begin{equation*}
(\mathbb{F}_q)^2 \subset \frac{\Delta(E)}{\Delta(E)}.
\end{equation*}
\end{theorem}
In this work, when $d=2$, we get the result of Theorem \ref{quotient} for all of $\mathbb{F}_q$ in the $q \equiv 3 \mod 4$ setting.

Recently, \cite{iosevich2023quotient} generalized the result of \cite{Iosevich_2019} for quadratic forms.  We state the result below after introducing the next definition.  For $r \in \mathbb{F}_q$ define
\begin{equation*}
V(r) := \# \left \{(a,b,c,d) \in E^4 : r = \frac{\|a-b\|^2}{\|c-d\|^2} \right \}. 
\end{equation*}
For simplicity, we state the relevant main results of \cite{iosevich2023quotient} in the specific case of the quadratic form $Q(x)=\|x\|^2$.
\begin{theorem}
\label{Firdav23}
(Theorem 1.6 in \cite{iosevich2023quotient})
Assume that $E \subset \mathbb{F}_q^d$, $d \geq 2$.
If $d \geq2$ is even and $|E|\geq4q^{d/2}$, then
\begin{equation*}
V(r) \geq \frac{5|E|^4}{48q} 
\end{equation*}
for all $r \in \mathbb{F}_q^*$.

If $d$ is odd and $|E|\geq \frac{11}{6}q^{(d+1)/2}$, then
\begin{align*}
V(r) \geq \frac{2|E|^4}{363q}.
\end{align*}
for all $r \in \mathbb{F}_q^*$
\end{theorem}
In Theorem \ref{Firdav23}, similar results to Theorem \cite{Iosevich_2019} are obtained with the addition of a count, $V(r)$.  In addition, like Theorem \cite{Iosevich_2019}, the proof of Theorem \ref{Firdav23} is quite sophisticated. 

The construction of the matrices $A_even$ and $A_odd$ from the proof of Theorem \ref{mainResult} allows us to improve the main result of \cite{bhowmik2023near}, which we state below.  Specifically, we allow for any $r \in \mathbb{F}_q^d$ instead of just $r \in (\mathbb{F}_q^{2})$.
\begin{theorem}
(Theorem 1.3 \cite{bhowmik2023near})
Suppose $r \in (\mathbb{F}_q^{2}) \setminus \{0\}$ and $\emptyset \not = A \subset \{
(i,j):1 \leq i < j \leq k+1\}$ where $k \geq 1$.  If $E \subset \mathbb{F}_q^d$ with $|E|\geq 2 kq^{d/2}$. then there exist $(x_1,...,x_{k+1}),(y_1,...,y_{k+1}) \in E^{k+1}$ such that $\|y_i -y_j\|^2=
r\|x_i-x_j\|^2$ if $(i,j) \in A$ and $x_i \not = x_j,y_i \not = y_j$ if $1 \leq i < j \leq k+1$.
\end{theorem}

\section{Results}
We use the group action approach and remove the restriction of $r \in (\mathbb{F}_q)^2$ from Theorem \ref{quotient}, allowing for any $r \in \mathbb{F}_q$.

\begin{theorem}
\label{mainResult}
Let $q$ be a prime power and $\E \subset \mathbb{F}_q^d$.  If $q \equiv 3 \mod 4$, then for each $r \in \mathbb{F}_q^*$,
\begin{equation*}
V(r)\geq \frac{|E|^4}{4q},
\end{equation*}
if $|E|\geq \sqrt{2}q$. 
\end{theorem}
The results of Theorem \ref{mainResult} are sharp.  For example, in section 1.1 of \cite{Iosevich_2019}, if $q=p^2$, in the $d=2$ setting, the set 
\begin{align*}
E= \mathbb{F}_p^2
\end{align*}
has size $q$, but
\begin{align*}
|\Delta(E)|=p,
\end{align*}
not $p^2$.

Using the proof techniques of Theorem \ref{mainResult}, we improve the main result of \cite{bhowmik2023near}, by allowing for any $r \in \mathbb{F}_q^d$ instead of just $r \in (\mathbb{F}_q^{2})$.
\begin{theorem}
\label{bhowmikResult}
Suppose $r \in (\mathbb{F}_q) \setminus \{0\}$ and $\emptyset \not = A \subset \{
(i,j):1 \leq i < j \leq k+1\}$ where $k \geq 1$.  If $E \subset \mathbb{F}_q^d$ with $|E|\geq 2 kq^{d/2}$, if $d$ is even or $|E|\geq 2 kq^{(d+1)/2}$, if $d$ is odd, then there exist $(x_1,...,x_{k+1}),(y_1,...,y_{k+1}) \in E^{k+1}$ such that $\|y_i -y_j\|^2=
r\|x_i-x_j\|^2$ if $(i,j) \in A$ and $x_i \not = x_j,y_i \not = y_j$ if $1 \leq i < j \leq k+1$.
\end{theorem}
Note the results are of Theorem \ref {bhowmikResult} are sharp.  In the case where $q=p^{2l}$, $p \equiv 3  \mod 4$, and $l \equiv 1 \mod 2$ the exponent $d/2$ is sharp (See section 6, page 9, \cite{bhowmik2023near}).

\begin{proof}
Proof of Theorem \ref{bhowmikResult}.  The proof follows almost verbatim of the combinatorial proof of Theorem 1.3 of \cite{bhowmik2023near} in the $r \not = 0,1$ case.  We assume that $r \not = 1$ is in $\mathbb{F}_q^*$.  Recall the matrix, $A$, from the proof Theorem \ref{mainResult}. In the proof of Theorem 1.3 of \ref{bhowmikResult}, page 8, we replace the set $tE$ with 
\begin{equation*}
AE = \{Av : v \in E\}
\end{equation*}
where $A = \sqrt{r} A_{odd}$, if $d$ is odd, and $A=\sqrt{r} A_{even}$, if $d$ is even, which we define below.  

\begin{equation}
A_{\text{even}}=\frac{1}{\sqrt{r}}\begin{bmatrix}
a  & -b  & &  & & & & &\\
b & a & &  & & & & &\\
   &   & a  & -b & &  & & & \\
   &   & b & a & &  & & &  \\
   &   &  &  & \ddots &  & &   \\
   &   &  &  &   & a  & -b & & &   \\
   &  & &  &   & b & a & & & \\
\end{bmatrix}.
\end{equation}
Since $d$ is even, we can choose $A$ to have $\frac{d}{2}$ blocks.  Note that $A_{even}$ is orthogonal by construction.  Denote $O(d)$, the group of orthogonal matrices in $\mathbb{F}_q^d$.
If $d\geq 3$ is odd we use the block matrix
\begin{equation*}
A_{\text{odd}}=
\begin{bmatrix}
A_{\text{even}} & 0  \\
0 & 1 \\
\end{bmatrix}.
\end{equation*}
noting $A_{odd}$ is orthogonal by construction.

As in the proof of Theorem 1.3 of \ref{bhowmikResult}, we define
\begin{align*}
H=\{(x,a): x \in AE \cap (E+a), a\in \mathbb{F}_q^d\},
\end{align*}
if $d$ is even.  If $d$ is odd we set
\begin{align*}
H=\{(x,a): x \in AE \cap (E+a), a\in \mathbb{F}_q^{d-1}\times \mathbb{F}_{q^2}\}.
\end{align*}
If $d$ is even,
then
\begin{align*}
|AE \cap (E+a)|\geq \frac{|E|^2}{q^d}
\end{align*}
by the averaging argument of Theorem 1.3 of \ref{bhowmikResult}.  Similarly if $d$ is odd,
\begin{align*}
|AE \cap (E+a)|\geq \frac{|E|^2}{q^{d+1}}
\end{align*}

We get $\{x_i\}_{i=1}^{k+2} \subset AE \cap (E+a)$ for some $a \in \mathbb{F}_q^{d-1}\times \mathbb{F}_{q^2}$ in the odd $d$ case and $a \in \mathbb{F}_q^d$ in the even $d$ case.  Thus \begin{align*}
y_i = x_i - a
\end{align*}
and 
\begin{align*}
A z_i = x_i
\end{align*}
for all $i=0,...,k+2$.

Thus 
\begin{align*}
\|y_i - y_j\|^2 = \|A(z_i-z_j\|^2
= r \|\frac{1}{\sqrt{r}}A(z_i-z_j) \|^2
=r\| (z_i-z_j)\|^2
\end{align*}
where the last equality above comes from the fact that $\sqrt{r}A$ is orthogonal.

\end{proof}

\begin{proof}
Proof of Theorem \ref{mainResult}.
The proof is similar to Theorem \ref{quotient}.  Assume that $r \in \mathbb{F}_q^*$.  Then $\sqrt{r} \in \mathbb{F}_{q^2}$.  For any characteristic, $p$, it is a basic fact from the theory of finite fields that there is some $a,b \in \mathbb{F}_q$ such that $a^2+b^2=r$.  We construct the diagonal matrix,
\begin{equation}
A=\frac{1}{\sqrt{r}}\begin{bmatrix}
a  & -b \\
b & a \\
\end{bmatrix}.
\end{equation}
Note that $A$ is orthogonal by construction.  Denote $O(2)$, the group of orthogonal matrices in $\mathbb{F}_q^2$.

Let $z \in \mathbb{F}_q^2$ and $\theta \in O(2)$.  Define the set,
\begin{equation*}
\eta_{\theta}(z)=
\{(u,v) \in E^2 :
u-\sqrt{r} \theta A v=z\}.
\end{equation*}
Consider the sum,
\begin{equation*}
\sum_{\theta \in O(d), z \in \mathbb{F}_q^d} \eta_{\theta}(z)= |E|^2|O(d)|.
\end{equation*}
 By the Cauchy Schwarz inequality,
\begin{equation*}
|E|^4|O(2)|^2= \left (\sum_{\theta \in O(2), z \in \mathbb{F}_q^2} \eta_{\theta}(z) \right)^2 \leq
\sum_{\theta,z}\eta_{\theta}(z)^2q^2|O(2)|.
\end{equation*}
It follows that,
\begin{align*}
\frac{|E|^4|O(2)|}{q^2}
\leq \sum_{\theta,z}\eta_{\theta}(z)^2.
\end{align*}

$\sum_{\theta,z}\eta_{\theta}(z)^2$  counts the number of  $(u,v,w,x) \in E^4$ such that $\|u-v\|^2
=r\|w-x\|^2$ modulo double counting by the stabilizer of $c-d$ for each $(c,d)\in E^2$.  The stabilizer of each $(c,d)\in E^2 \backslash \text{diag}(E)$ has $|O(1)|$ elements (see page 3 of \cite{bennett2013group}).

The number of tuples, $(u,v,w,x)$ in the sum, $\sum_{\theta,z}\eta_{\theta}(z)^2$, where $u = v$ and $w=x$ are counted $\sum_{\theta \in O(2), z \in \mathbb{F}_q^2} \eta_{\theta}(z)= |E|^2|O(2)|$ times.  

Define 
\begin{equation*}
N_0 = \#\{(u,v,w,x) \in E^4 : \|u-v\|^2
= r\|w-x\|^2=0 \text{ and }u-v,w-x \not = 0\}.
\end{equation*}

Since $q \equiv 3 \mod 4$, $N_0=0$. Thus, if $|E|\geq \sqrt{2}q$, 
\begin{equation*}
V(r) \geq \frac{1}{|O(1)|} \left (\sum_{\theta,z}\eta_{\theta}(z)^2
-|E|^2|O(2)| - |O(1)|N_0 \right )
\end{equation*}
\begin{equation*}
=  \frac{1}{2}\left (\sum_{\theta,z}\eta_{\theta}(z)^2
-|E|^2|O(2)| \right)
\geq  \frac{1}{2} \left (\frac{|E|^4(q+1)}{q^2} - |E^2|(q+1) \right ).
\end{equation*}
From above, if $|E| \geq \sqrt{2}q$, then $V(r)\geq \frac{|E|^4(q+1)}{4q^2}$, and we're done.
\end{proof}

\bibliography{iclr2024_conference}
\bibliographystyle{iclr2024_conference}

\appendix

\end{document}